\documentclass{etds}
\usepackage{amsmath,amssymb,bm,color}

\newtheorem{theo}{Theorem}

\newtheorem{lemma}{Lemma}

\begin{document}
\ETDS{0}{1}{1}{2011}

\runningheads{I. Assani and K. Presser}{Pointwise Characteristic Factors}

\title{Pointwise characteristic factors for the multiterm return times theorem}

\author{Idris Assani\affil{1}, Kimberly Presser\affil{2}}

\address{\affilnum{1}\ Department of Mathematics,
UNC Chapel Hill, NC 27599, \\ \email{assani@math.unc.edu} \\ \affilnum{2}\ Department of Mathematics, Shippensburg University, PA 27517 \\ \email{kjpres@ship.edu}}

\recd{December 30, 2010}

\begin{abstract}
This paper is an update and extension of a result the authors first proved in 2003.  The goal of this paper is to study factors which are known to be $L^2$-characteristic for certain nonconventional averages and prove that these factors are pointwise characteristic for the multiterm return times averages.

\begin{center}In memory of Dan Rudolph.\end{center}
\end{abstract}

\section{Introduction}
A major result in ergodic theory in the late 1980's was the proof of the return times theorem by J. Bourgain \cite{RTT1} (which was later simplified by J. Bourgain, H. Furstenberg, Y. Katznelson, D. Ornstein in \cite{RTT2}).  This theorem created a key strengthening of the Birkhoff's Pointwise Ergodic Theorem \cite{Birk}.

\begin{theo}\label{RTT1}
Let $(X,\mathcal{F},\mu,T)$ be an ergodic dynamical system of finite measure and $f \in L^{\infty}(\mu)$.  Then there exists a set $X_f \subset X$ of full measure such that for any other ergodic dynamical system $(Y,\mathcal{G}, \nu, S)$ with $\nu(Y)<\infty$ and any $g \in L^{\infty}(\nu)$:
$$\frac{1}{N}\sum_{n=1}^Nf(T^nx)g(S^ny)$$
converges $\nu$-a.e. for all $x \in X_f$.
\end{theo}

Note that the set $X_f$ depends not only on the function $f$ chosen, but on the transformation $T$ in our dynamical system.  In the BFKO proof \cite{RTT2} of the return times theorem, one of the keys to the argument was to decompose the given function using the Kronecker factor in order to prove the result independently for both the eigenfunctions and those functions in the orthocomplement of the Kronecker factor.

Using factors in convergence proofs in ergodic theory has long been a very useful tool.  The notion of a characteristic factor is originally due to H. Furstenberg and is explicitly defined by H. Furstenberg and B. Weiss in \cite{FW}.

\proc{Definition:}
When the limiting behavior of a non-conventional ergodic average for $(X,\mathcal{F}, \mu, T)$ can be reduced to that of a factor system $(Y, \mathcal{G}, \nu, T)$, we shall say that the latter is a \textbf{characteristic factor} of the former.
\medbreak

For each type of average under consideration, one will have to specify what is meant by reduced in the given case.  In the case of H. Furstenberg and B. Weiss \cite{FW}, they define the notion of characteristic factor for averages of the type
$$\frac{1}{N}\sum_{n=1}^N \left(f\circ T^n\right)\left(g \circ T^{n^2}\right).$$
Therefore their specific definition of characteristic factor is as follows.

\proc{Definition:}
If $\{p_1(n), p_2(n), \ldots, p_k(n)\}$ are $k$ integer-valued sequences, and $(Y,\mathcal{G}, \nu, T)$ is a factor of a system $(X, \mathcal{F}, \mu, T)$, we say that $\mathcal{G}$ is a \textbf{characteristic factor for the scheme} $\{p_1(n), p_2(n), \ldots, p_k(n)\}$, if for any $f_1, f_2, \ldots, f_k \in L^{\infty}(\mu)$ we have
$$\frac{1}{N} \sum_{n=1}^N \left[f_1 \circ T^{p_1(n)}\cdots f_k \circ T^{p_k(n)}-\mathbb{E}(f_1 | \mathcal{G}) \circ T^{p_1(n)}\cdots \mathbb{E}(f_k|\mathcal{G})\circ T^{p_k(n)}\right]$$
converges to $0$ in $L^2(\mu)$.
\medbreak

In 1998, D. Rudolph \cite{RudolphFGS} extended the return times theorem to averages with more than two terms with his proof of the multiterm return times theorem.  His proof answered one of the questions  on the return times raised by I. Assani\footnote{These questions were brought up during D. Rudolph's visit to UNC-CH in 1991 while he was working on his joinings proof of Bourgain's return times theorem \cite{RudolphAJP}.} who proved the same result for weakly mixing systems in \cite{AssaniMRA}.

\begin{theo}[Multiterm Return Times Theorem]\label{MRTT}
Let $k$ be any positive integer.  For any ergodic dynamical system $(X,\mathcal{F}, T, \mu)$ and any $f \in L^{\infty}(\mu)$, there exists a set of full measure $X_{f}$ in $X$ such that if $x \in X_{f}$ for any other dynamical system $(Y_1,\mathcal{G}_1, S_1,\nu_1)$ and any $g_1 \in L^{\infty}(\nu_1)$ there exists a set of full measure $Y_{g_1}$ in $Y_1$ such that if $y_1 \in Y_{g_1}$then \ldots for any other dynamical system $(Y_{k-1},\mathcal{G}_{k-1}, S_{k-1}, \nu_{k-1})$ and any $g_{k-1} \in L^{\infty}(\nu_{k-1})$ there exists a set of full measure $Y_{g_{k-1}}$ in $Y_{k-1}$ such that if $y_{k-1} \in Y_{g_{k-1}}$ for any other dynamical system $(Y_{k},\mathcal{G}_{k}, S_{k},\nu_{k})$ the average:
$$\frac{1}{N}\sum_{n=1}^N f(T^{n}x)g_{1}(S_{1}^{n}y_{1})g_{2}(S_{2}^{n}y_{2})\cdots g_{k}(S_{k}^{n}y_{k})$$
converges $\nu_{k}$-a.e..
\end{theo}

As we are interested in finding characteristic factors for ergodic dynamical systems this theorem is quoted here with the extra assumption of ergodicity for the dynamical system $(X,\mathcal{F}, T, \mu)$.  D. Rudolph's proof of the multiterm return times theorem utilized the method of joinings and fully generic sequences.  This led to an elegant proof of the theorem which avoided the study of the factor of the $\sigma$-algebra which was characteristic for the averages.  So the higher order version of the Kronecker factor $\mathcal{K}$ which had been key to the BFKO \cite{RTT2} proof was not needed in D. Rudolph's argument.  This paper seeks to determine what factors serve a role similar to the Kronecker factor $\mathcal{K}$ in this multiterm setting.

For our purposes we define the notion of pointwise characteristic factors for the multiterm return times averages as follows.

\proc{Definition:}
Consider $(X, \mathcal{F}, \mu, T)$ a measure preserving system.  The factor $\mathcal{A}$ is \textbf{pointwise characteristic for the $k$-th return times averages} if for each $f \in L^{\infty}(\mu)$ we can find a set of full measure $X_f$ such that for each $x \in X_f$, for any other dynamical system $(Y_1,\mathcal{G}_1, S_1, \nu_1)$ and any $g_1 \in L^{\infty}(\nu_1)$, there exists a set of full measure $Y_{g_1}$ such that for each $y_1$ in $Y_{g_1}$ then $\ldots$ for any other dynamical system $(Y_{k-1},\mathcal{G}_{k-1}, S_{k-1}, \nu_{k-1})$ and any $g_{k-1} \in L^{\infty}(\nu_{k-1})$, there exist a set of full measure $Y_{g_{k-1}}$ in $Y_{k-1}$ such that if $y_{k-1} \in Y_{g_{k-1}}$ for any other dynamical system $(Y_k,\mathcal{G}_k, S_k, \nu_k)$ for $\nu_k$-a.e. $y_k$ the average
$$ \frac{1}{N}\sum_{n=1}^N \left[f(T^{n}x)-\mathbb{E}(f | \mathcal{A})(T^nx)\right]g_{1}(S_{1}^{n}y_{1})g_{2}(S_{2}^{n}y_{2})\cdots g_{k}(S_{k}^{n}y_{k})$$
converges to 0.
\medbreak

In looking for potential characteristic factors for the general multiterm return times averages we consider the factors first used by H. Furstenberg to prove Szemer\'{e}di's Theorem \cite{FurstEBO}.  These factors are called $k$ step distal factors in \cite{FurstEBO}. We denote these factors (which will be further defined in Section 2) as $\mathcal{A}_k$ using the notation from \cite{AssaniCFF} where these factors were shown to be $L^2$-characteristic for the averages
$$\frac{1}{N}\sum_{n=1}^N \prod_{i=1}^I f_i\circ T^{in}.$$

While the norm convergence of averages for $L^2$-characteristic factors can sometimes lead to pointwise characteristic properties, this is not always guaranteed to be the case.  Thus it is of consequence to look at pointwise convergence in addition to investigating factors with respect to the norm convergence.

We will show that these $\mathcal{A}_k$ factors can be characterized in an inductive way by seminorms using Lemma 1.3 of \cite{RudolphEOT}.  Using these seminorms we will prove our first result.\footnote{This approach was used in two 2003 unpublished papers of the first author (\cite{AssaniCFF} and what was ultimately combined into the published paper \cite{AssaniPCAC}).  The first author thanks C. Demeter and N. Frantzikinakis for pointing out to him that the factors he defined with these seminorms were in fact the ones introduced by H. Furstenberg in \cite{FurstEBO}. A careful look at Theorem 10.2 in \cite{FurstEBO} indicates that the k step distal factors are $L^2$ characteristic for the Furstenberg averages.}

\begin{theo}\label{Main}
The factors $\mathcal{A}_k$ are pointwise characteristic for the multiterm return times averages.  More precisely, using the seminorms defining the $\mathcal{A}_k$ we can find pointwise uniform upper bounds of the multiterm return times averages.
\end{theo}

The study of the nonconventional Furstenberg averages has seen important progress being made in the last seven years. In \cite{HK-NEA} and \cite{Ziegler} the Host-Kra-Ziegler factors $\mathcal{Z}_k$ were created independently by B. Host, B. Kra and T. Ziegler and were shown to be characteristic in $L^2$ norm for the Furstenberg averages.  Using these factors we prove our second result.

\begin{theo}\label{HKZ}
Let $(X, \mathcal{F}, \mu, T)$ be an ergodic measure preserving system.  The Host-Kra-Ziegler factors $\mathcal{Z}_k$ are pointwise characteristic for the multiterm return times averages.
\end{theo}

As the $\mathcal{Z}_k$ factors are smaller than the factors $\mathcal{A}_k$, and thus $\mathcal{A}_k^{\perp} \subseteq \mathcal{Z}_k^{\perp}$, the fact that the $\mathcal{A}_k$ factors are pointwise characteristic for the multiterm return times averages is a consequence of Theorem \ref{HKZ}. But in our proof of Theorem \ref{Main} using the seminorm defining the factors $\mathcal{A}_k$ we obtain pointwise uniform upper bounds of the multiterm return times averages. With the $\mathcal{Z}_k$ factors we do not have such pointwise estimates. The uniform upper bounds are derived after integration combined with a $\limsup$ argument.

An unusual feature of our proof is that we use the previously established convergence result of D. Rudolph and use that to prove that the factors we are studying are characteristic.

\section{The $k$-step distal factors $\mathcal{A}_k$ factors are pointwise characteristic for the multiterm return times averages}

Let $(X, \mathcal{F}, \mu, T)$ be an ergodic dynamical system on a probability measure space.  The factors $\mathcal{A}_k$ are defined in the following inductive way.

\proc{Definition:}
\begin{itemize}
\item The factor $\mathcal{A}_0$ is equal to the trivial $\sigma$-algebra $\{X,\emptyset\}$
\item For $k\geq 0$ the factor $\mathcal{A}_{k+1}$ is characterized by the following.  A function $f\in \mathcal{A}_{k+1}^{\perp}$ if and only if
$$N_{k+1}(f)^4 := \lim_H\frac{1}{H}\sum_{h=1}^H \left\|\mathbb{E}(f\cdot f\circ T^h | \mathcal{A}_{k})\right\|_2^2 = 0$$
\end{itemize}
\medbreak

Note that the factor $\mathcal{A}_1$ is the Kronecker factor of our ergodic transformation $T$ because
$$N_1(f)^4 =\lim_H \frac{1}{H} \sum_{h=1}^H \left\|\mathbb{E}(f\cdot f\circ T^h | \mathcal{A}_0)\right\|_2^2 = \lim_H \frac{1}{H}\sum_{h=1}^H \left|\int f\cdot f \circ T^h d\mu\right|^2.$$
In Lemma \ref{seminorms} we will verify that the $\mathcal{A}_k$ as defined above do actually form well-defined factors.

We will want to verify that these $\mathcal{A}_k$ are maximal isometric extensions.  There are several equivalent ways of expressing this feature.  We will be using the terminology found on pages 373-374 of \cite{RudolphEOT} to specify how these factors form maximal isometric extensions.  Furstenberg has shown in \cite{FurstDistal} that for any ergodic dynamical system $(X, \mathcal{F}, \mu, T)$ and any $T$-invariant factor $\mathcal{B}$ there is a unique maximal factor action $\mathcal{KB} \subseteq \mathcal{F}$ which contains $\mathcal{B}$ so that in the Rohlin representation of $(T,\mathcal{KB})$, the space $\mathcal{Z}_2$ can be taken to be a compact metric space with isometric factor maps.  This factor $\mathcal{KB}$ arises from the invariant algebras of the relatively independent joinings of $(X,\mathcal{F},\mu,T)$ over the factor $\mathcal{B}$.  This is precisely the notion of maximal isometric extension referred to in the next lemma.  

\begin{lemma}\label{seminorms}
Let $(X, \mathcal{F}, \mu, T)$ be an ergodic dynamical system on a probability measure space.  For $k\geq 2$ for each function $f\in L^{\infty}(\mu)$  the quantities $N_k(f)$ are well defined. Furthermore, they characterize factors of $T$ which are successive maximal isometric extensions.
\end{lemma}

\proc{Proof:}
Let us consider a general factor $\mathcal{A}$ of $T$ and $\mathbb{E}(\cdot|\mathcal{A})$ the projection onto this factor. The relatively independent joining of $T\times T$ over the factor $\mathcal{A}$ is the measure $\mu_{\mathcal{A}}$ defined for $f,$ $g$ bounded functions as
$$\int f\times g d\mu_{\mathcal{A}}:=\int \mathbb{E}(f|\mathcal{A})\mathbb{E}(g|\mathcal{A})d\mu.$$
By Birkhoff's ergodic theorem applied to $T\times T$ and the invariant measure $\mu_{\mathcal{A}}$ we have
\begin{eqnarray*}
\lim_H \frac{1}{H}\sum_{h=1}^H \left\|\mathbb{E}(f\cdot f\circ T^h|\mathcal{A})\right\|_2^2 & = & \lim_H\frac{1}{H}\sum_{h=1}^H \int(f \cdot f\circ T^h)(x) (f\cdot f\circ T^h)(y)d\mu_{\mathcal{A}} \\
& = & \left\|\mathbb{E}(f\times f|\mathcal{I}_{\mathcal{A}})\right\|_{L^2(\mu_{\mathcal{A}})}^2.
\end{eqnarray*}
where $\mathcal{I}_{\mathcal{A}}$ is the $T\times T$-$\mu_{\mathcal{A}}$ invariant $\sigma$-algebra.

If we denote by $N(f)$ the quantity
$$N(f)^4 = \lim_H\frac{1}{H}\sum_{h=1}^H \left\|\mathbb{E}(f \cdot f\circ T^h|\mathcal{A})\right\|_2^2$$
then Lemma 1.3 in \cite{RudolphEOT} tells us that $N(f)=0$ if and only if $\mathbb{E}(f|\mathcal{KA})=0$ where $\mathcal{KA}$ is the maximal isometric extension of $\mathcal{A}.$

Using these observations one can characterize the successive maximal isometric extensions.  The trivial $\sigma$-algebra is $\mathcal{A}_0$. Then we define $\mathcal{A}_1 = \mathcal{KA}_0$, $\mathcal{A}_2 = \mathcal{KA}_1$ and more generally $\mathcal{A}_{k+1} = \mathcal{KA}_{k}.$ The seminorms characterizing these factors are well defined as $N_k(f)$ where
$$N_k(f)^4 = \lim_H \frac{1}{H}\sum_{h=1}^H \left\|\mathbb{E}(f\cdot f\circ T^h|\mathcal{A}_k)\right\|_2^2.$$
\ep
\medbreak

In order to simplify the inductive parts of our argument, we first clarify the techniques that we will use in a series of small lemmas.  This next lemma relies on an application of the spectral theorem which allows us to alternate between Wiener-Wintner and return times averages in our inductive argument.

\begin{lemma}\label{WWtype}
Let $\{a_n\}$ be a sequence of complex numbers. If
$$\sup_N\frac{1}{N}\sum_{n=1}^N|a_n|^2<\infty \textrm{ and } \sup_{\epsilon}\left|\frac{1}{N}\sum_{1}^{N}a_ne^{2\pi in\epsilon}\right|\rightarrow 0,$$
then
$$\frac{1}{N}\sum_{1}^{N}a_ng(S^ny)\rightarrow 0$$
in $L^2(\nu)$ for all measure-preserving systems $(Y,\mathcal{G},S,\nu)$.
\end{lemma}

\proc{Proof:}
This follows immediately from the proof of Theorem 3.1 in \cite{AssaniWWE}.
\ep
\medbreak

Next, we will use the following lemma which is an easy consequence of the Van der Corput lemma \cite{KN}.  It will help us simplify the Wiener-Wintner averages which will appear in the inductive argument.

\begin{lemma}\label{VDCtype}
There exists an absolute constant $C$ such that for any $\{a_n\}$ bounded sequence of complex numbers and any positive integer $N$ we have
$$\sup_{\epsilon}\left|\frac{1}{N}\sum_{n=1}^{N}a_ne^{2\pi in\epsilon}\right|^2 \leq C \left(\frac{1}{H}+\frac{1}{H}\sum_{h=1}^{H}\left|\frac{1}{N}\sum_{n=1}^{N-h}a_n\overline{a_{n+h}}\right|\right)$$
for $1 \leq H \leq N$.
\end{lemma}

The following lemma will be useful in establishing the basis step of our next theorem. It gives a pointwise upper bound for the return times averages for two terms, case studied in the BFKO \cite{RTT2} paper.

\begin{lemma}\label{Upbound0}
Let $(X, \mathcal{F}, \mu, T)$ be an ergodic dynamical system and $f\in L^{\infty}(\mu).$ Let us denote by $\mathcal{K}_T$ its Kronecker factor. Then there exists a universal set $X_f$ depending only on $f$ and the system $(X,\mathcal{F}, \mu, T)$ such that for any dynamical system $(Y_1,\mathcal{G}_1, S_1, \nu_1)$ and any $g_1\in L^{\infty}(\nu_1)$ we have
\begin{equation}
 \limsup_N \left|\frac{1}{N}\sum_{n=1}^N f(T^nx)g_1(S_1^ny)\right|\leq \|\mathbb{E}(f| \mathcal{K}_T)\|_2\|g_1\|_{\infty}.
\end{equation}
\end{lemma}

\proc{Proof:}\
By the BFKO return times theorem, we know that the Kronecker factor is pointwise characteristic.  So we have on a universal set $X_f$ of full measure
$$\limsup_N \left|\frac{1}{N}\sum_{n=1}^N f(T^nx)g_1(S_1^ny_1)\right|= \limsup_N \left|\frac{1}{N}\sum_{n=1}^N \mathbb{E}(f| \mathcal{K}_T)(T^nx)g_1(S_1^ny_1)\right|,$$ $\nu_1$ a.e.,
for any dynamical system $(Y_1,\mathcal{G}_1, S_1, \nu_1)$ and any $g_1\in L^{\infty}(\nu_1)$.  

Therefore for any dynamical system $(Y_1,\mathcal{G}_1, S_1, \nu_1)$ and any $g_1\in L^{\infty}(\nu_1)$,
$$\limsup_N \left|\frac{1}{N}\sum_{n=1}^N f(T^nx)g_1(S_1^ny_1)\right| \leq \limsup_N \left|\frac{1}{N}\sum_{n=1}^N \mathbb{E}(f| \mathcal{K}_T)(T^nx)\right|\|g_1\|_{\infty}.$$
Using the Cauchy Schwartz inequality we have then
$$\limsup_N \left|\frac{1}{N}\sum_{n=1}^N f(T^nx)g_1(S_1^ny_1)\right|\leq \limsup_N \left(\frac{1}{N}\sum_{n=1}^N |(\mathbb{E}(f|\mathcal{K}_T)(T^nx)|^2\right)^{1/2}\|g_1\|_{\infty}.$$

Then by applying Birkhoff's pointwise ergodic theorem we have 
$$\limsup_N \left|\frac{1}{N}\sum_{n=1}^N f(T^nx)g_1(S_1^ny_1)\right|\leq \left(\int |(\mathbb{E}(f|\mathcal{K}_T)|^2 d\mu\right)^{1/2}\|g_1\|_{\infty},$$
which is the upper bound announced in this lemma.
\ep
\medbreak

We will prove our first main result, Theorem \ref{Main}, in the course of proving the following more detailed statement.

\begin{theo}\label{A_k}
Let $k$ be any positive integer.  For any ergodic dynamical system $(X,\mathcal{F}, \mu, T)$ and for each $f \in L^{\infty}(\mu)$ we can find a set of full measure $X_f$ such that for each $x \in X_f$, for any other dynamical system $(Y_1,\mathcal{G}_1, S_1, \nu_1)$ and any $g_1 \in L^{\infty}(\nu_1)$ with $\|g_1\|_{\infty} \leq 1$, there exists a set of full measure $Y_{g_1}$ such that for each $y_1$ in $Y_{g_1}$ then $\ldots$ for any other dynamical system $(Y_{k-1},\mathcal{G}_{k-1}, S_{k-1}, \nu_{k-1})$ and any $g_{k-1} \in L^{\infty}(\nu_{k-1})$ with $\|g_k\|_{\infty} \leq 1$ there exist a set of full measure $Y_{g_{k-1}}$ in $Y_{k-1}$ such that if $y_{k-1} \in Y_{g_{k-1}}$ for any other dynamical system $(Y_k,\mathcal{G}_k, S_k, \nu_k)$ for $\nu_k$-a.e. $y_k$
\begin{itemize}
    \item the average
        \begin{equation}\label{RTT}
            \frac{1}{N}\sum_{n=1}^N \left[f(T^{n}x)-\mathbb{E}(f | \mathcal{A}_k)(T^nx)\right]g_{1}(S_{1}^{n}y_{1})g_{2}(S_{2}^{n}y_{2})\cdots g_{k}(S_{k}^{n}y_{k})
        \end{equation}
        converges to 0.
    \item Thus for $f\in \mathcal{A}_k^{\perp}$ the average
$$\frac{1}{N}\sum_{n=1}^N f(T^{n}x)g_{1}(S_{1}^{n}y_{1})g_{2}(S_{2}^{n}y_{2})\cdots g_{k}(S_{k}^{n}y_{k})$$
converges to 0 $\nu_k$-a.e..
    \item Also we have the following pointwise upper bound for our limit
        \begin{equation}\label{Upbound}
            \limsup_N\left|\frac{1}{N}\sum_{n=1}^N f(T^nx)g_1(S_1^ny)g_{2}(S_{2}^{n}y_{2})\cdots g_{k}(S_{k}^{n}y_{k})\right|^2\leq CN_{k+1}(f)^2
        \end{equation}
\end{itemize}
\end{theo}

\proc{Proof:}
The basis step for the induction of the statement in (\ref{RTT}) was done in the BFKO \cite{RTT2} proof of Bourgain's Return Times Theorem . Here it was shown that $\mathcal{A}_1=\mathcal{K}_T$ was pointwise characteristic for averages of the type
\begin{equation}\label{Average}
\frac{1}{N}\sum_{n=1}^N f(T^nx)g_1(S_1^ny_1).
\end{equation}

In Lemma \ref{Upbound0}, we showed that the quantity $\|\mathbb{E}(f|\mathcal{A}_1)\|_2\|g_1\|_{\infty}$ is a pointwise upper bound for the $\limsup$ of the absolute value of the averages where $\mathcal{K}_{S_1}$ is the Kronecker factor for $S_1.$  This last term is itself less than $\|\mathbb{E}(f|\mathcal{A}_1)\|_2$ because  $\|g_1\|_{\infty}\leq 1.$  Thus we have reached the inequality
\begin{equation}\label{Bound}
 \limsup_N\left|\frac{1}{N}\sum_{n=1}^N f(T^nx)g_1(S_1^ny)\right|^2\leq \|\mathbb{E}(f|\mathcal{A}_1)\|_2^2.
\end{equation}

We want to get a better upper bound namely $CN_2(f)$ where $C$ is an absolute constant. To this end we apply the Van der Corput lemma to obtain
$$\limsup_N \left|\frac{1}{N} \sum_{n=1}^N f(T^nx) g_1(S_1^ny)\right|^2 \leq$$
$$\limsup_N C\left[\frac{1}{H} + \left(\frac{1}{H}\sum_{h=1}^H \left|\frac{1}{N}\sum_{n=1}^{N-h}f(T^nx)f(T^{n+h}x)g_1(S_1^ny)g_1(S_1^{n+h}y_1)\right|\right)\right] \leq$$
$$C\left[\frac{1}{H} + \left(\frac{1}{H}\sum_{h=1}^H \limsup_N\left|\frac{1}{N}\sum_{n=1}^{N-h}f(T^nx)f(T^{n+h}x)g_1(S_1^ny)g_1(S_1^{n+h}y_1)\right|\right)\right].$$

Applying the inequality (\ref{Bound}) to each of the functions $f\cdot f\circ T^h$ yields
$$\limsup_N \left|\frac{1}{N} \sum_{n=1}^N f(T^nx) g(S^ny)\right|^2 \leq  C\left[\frac{1}{H} + \left(\frac{1}{H}\sum_{h=1}^H \|\mathbb{E}(f\cdot f\circ T^h|\mathcal{A}_1)\|_2^2\right)^{1/2}\right].$$

By taking the limit with $H$ we get the better estimate
$$\limsup_N \left| \frac{1}{N} \sum_{n=1}^N f(T^nx) g(S^ny)\right|^2 \leq C N_2(f)^2.$$
which shows clearly that $\mathcal{A}_1$ satisfies (\ref{Upbound}) in the basis step.

Assume that for any $f \in L^{\infty}(\mu)$ and $1 \leq j < k$ we can find sets $X_f$ of full measure such that if $x \in X_f$, then for any other dynamical system $(Y_1,\mathcal{G}_1, S_1, \nu_1)$ and any $g_1 \in L^{\infty}(\nu_1)$ with $\|g_1\|_{\infty} \leq 1$, there exists a set of full measure $Y_{g_1}$ such that for each $y_1$ in $Y_{g_1}$ then $\ldots$ for any other dynamical system $(Y_{j-1},\mathcal{G}_{j-1}, S_{j-1}, \nu_{j-1})$ and any $g_{j-1} \in L^{\infty}(\nu_{j-1})$ with $\|g_{j-1}\|_{\infty} \leq 1$ there exist a set of full measure $Y_{g_{j-1}}$ in $Y_{j-1}$ such that if $y_{j-1} \in Y_{g_{j-1}}$ for any other dynamical system $(Y_j,\mathcal{G}_j, S_j, \nu_j)$ and any $g_{j} \in L^{\infty}(\nu_{j})$ with $\|g_{j}\|_{\infty} \leq 1$ for $\nu_j$-a.e. $y_j$ we have
\begin{itemize}
    \item the average
        $$ \frac{1}{N}\sum_{n=1}^N \left[f(T^{n}x)-\mathbb{E}(f | \mathcal{A}_j)(T^nx)\right]g_{1}(S_{1}^{n}y_{1})\cdots g_{j}(S_{j}^{n}y_{j})$$
    converges to 0.
    \item  Also we have the upper bound
        $$\limsup_N\left|\frac{1}{N}\sum_{n=1}^N f(T^nx)g_1(S_1^ny)g_{2}(S_{2}^{n}y_{2})\cdots g_{j}(S_{j}^{n}y_{j})\right|^2\leq CN_{j+1}(f)^2$$
\end{itemize}

\begin{lemma}\label{BN}
Let $f$ be an element of $f \in L^{\infty}$ and let $g_i$, $S_i$ and $y_i$ be as defined in the preceding paragraph.  If
$$B_N=\sup_{\epsilon}\left|\frac{1}{N}\sum_{n=1}^{N}f(T^{n}x)g_{1}(S_{1}^{n}y_{1})\cdots g_{k-1}(S_{k-1}^{n}y_{k-1})e^{2\pi in\epsilon}\right|^2$$
then
$$\limsup_N B_N \leq CN_k(f)^2$$
for some absolute constant $C$.  Here the constant $C$ is independent of the $f$, $g_i$, $S_i$ and $y_i$.
\end{lemma}

\proc{Proof:}
By Lemma \ref{VDCtype}, there exists a constant $C$ such that for $1 \leq H\leq N$
\begin{eqnarray*}
B_N & \leq & C\Bigg(\frac{1}{H}+\frac{1}{H}\sum_{h=1}^{H}\Bigg|\frac{1}{N}\sum_{n=1}^{N-h}(f\cdot \overline{f \circ T^h})(T^{n}x)\\
& & \cdot(g_1\cdot \overline{g_1 \circ
S_1^h})(S_{1}^{n}y_{1})\cdots (g_{k-1}\cdot \overline{g_{k-1} \circ S_{k-1}^h})(S_{k-1}^{n}y_{k-1})\Bigg|\Bigg).
\end{eqnarray*}

From our inductive hypothesis, we know that for each $h$ there is a set of full measure $X_{f\cdot \overline{f \circ T^h}}$ on which
\begin{eqnarray*}
\Bigg|\frac{1}{N}\sum_{n=1}^{N-h}\left[(f\cdot \overline{f \circ T^h})(T^{n}x)-(\mathbb{E}(f\cdot \overline{f \circ
T^h}|\mathcal{A}_{k-1})(T^{n}x)\right] & &\\
\cdot(g_1\cdot \overline{g_1 \circ S_1^h})(S_{1}^{n}y_{1})\cdots (g_{k-1}\cdot \overline{g_{k-1} \circ S_{k-1}^h})(S_{k-1}^{n}y_{k-1})\Bigg| & \rightarrow & 0.
\end{eqnarray*}

Therefore, the intersection of these sets $X_{f\cdot \overline{f \circ T^h}}$ over $h$ gives a set of full measure $\widehat{X_f}$ on which
\begin{eqnarray*}
\limsup_N B_N & \leq & \limsup_NC\Bigg(\frac{1}{H}+\frac{1}{H}\sum_{h=1}^{H}\Bigg|\frac{1}{N}\sum_{n=1}^{N-h}(\mathbb{E}(f\cdot \overline{f \circ T^h}|\mathcal{A}_{k-1})(T^{n}x)\\
& & \cdot(g_1\cdot \overline{g_1 \circ S_1^h})(S_{1}^{n}y_{1})\cdots (g_{k-1}\cdot \overline{g_{k-1}
\circ S_{k-1}^h})(S_{k-1}^{n}y_{k-1})\Bigg|\Bigg)
\end{eqnarray*}
for all $H$.

The Cauchy-Schwartz inequality gives us
\begin{eqnarray*}
\limsup_NB_N & \leq & \limsup_N C \Bigg(\frac{1}{H}+\frac{1}{H}\sum_{h=1}^{H}\Bigg(\frac{1}{N}\sum_{n=1}^{N-h}\left|(\mathbb{E}(f\cdot \overline{f \circ T^h}|\mathcal{A}_{k-1})(T^{n}x)\right|^2\\
&  & \cdot\left|(g_1\cdot \overline{g_1 \circ S_1^h})(S_{1}^{n}y_{1})\right|^2 \cdots \left|(g_{k-1}\cdot
\overline{g_{k-1} \circ S_{k-1}^h})(S_{k-1}^{n}y_{k-1})\right|^2\Bigg)^{\frac{1}{2}}\Bigg) \\
& \leq & \limsup_N C\Bigg(\frac{1}{H}+\frac{\|g_{1}\|_{\infty}^2\ldots \|g_{k-1}\|_{\infty}^2}{H}\cdot \\
& & \sum_{h=1}^{H}\left(\frac{1}{N}\sum_{n=1}^{N-h}\left|(\mathbb{E}(f\cdot \overline{f \circ T^h}|\mathcal{A}_{k-1})(T^{n}x)\right|^2\right)^{\frac{1}{2}}\Bigg)
\end{eqnarray*}

By Birkhoff's Pointwise Ergodic Theorem we know that there is a set of full measure $X_{k-1}$ on which for each $h$ the average over $n$ in the above inequality converges to
$$\int\left|(\mathbb{E}(f\cdot \overline{f \circ T^h}|\mathcal{A}_{k-1})\right|^2d\mu=\left\|\mathbb{E}(f\cdot \overline{f \circ T^h}|\mathcal{A}_{k-1})\right\|^2_2.$$
Therefore on the set of full measure $X_f = \widehat{X_f} \bigcap X_{k-1}$
\begin{eqnarray*}
\limsup_NB_N & \leq & \frac{C}{H}+\frac{C}{H}\sum_{h=1}^{H}\left\|\mathbb{E}(f\cdot \overline{f \circ T^h}|\mathcal{A}_{k-1})\right\|_2 \\
& \leq & C\lim_H \left(\frac{1}{H}+\left(\frac{1}{H}\sum_{h=1}^{H}\left\|\mathbb{E}(f\cdot \overline{f \circ T^h}|\mathcal{A}_{k-1})\right\|_2^2\right)^{\frac{1}{2}}\right) \\
& = & C\cdot N_k(f)^2.
\end{eqnarray*}
\ep
\medbreak

As functions $f$ in $\mathcal{A}_k^{\perp}$ are characterized by the property that $N_k(f)=0$, Lemma \ref{BN} implies that when $f$ is an element of $L^{\infty}(\mu)\bigcap\mathcal{A}_k^{\perp}$ we have
$$\limsup_NB_N =0$$
on the set of full measure $X_f=\widehat{X_f}\bigcap X_{k-1}$.  Therefore
$$\sup_{\epsilon}\left|\frac{1}{N}\sum_{n=1}^{N}f(T^{n}x)g_{1}(S_{1}^{n}y_{1})\cdots g_{k-1}(S_{k-1}^{n}y_{k-1})e^{2\pi in\epsilon}\right|$$
converges to 0 $\mu$-a.e..  Hence by an application of Lemma \ref{WWtype}, we know that for any other dynamical system $(Y_k,\mathcal{G}_k,S_k,\nu_k)$ and any $g \in L^{\infty}(\nu_k)$
\begin{equation}\label{norm}
\frac{1}{N}\sum_{n=1}^N f(T^{n}x)g_{1}(S_{1}^{n}y_{1})\cdots g_{k}(S_{k}^{n}y_{k})
\end{equation}
converges to 0 in $L^2(\nu_k)$. As pointwise convergence of the average in Equation (\ref{norm}) follows from Theorem \ref{MRTT}, we have
$$\frac{1}{N}\sum_{n=1}^N f(T^{n}x)g_{1}(S_{1}^{n}y_{1})\cdots g_{k}(S_{k}^{n}y_{k})$$
converges to 0 $\nu_k$-a.e., when $f$ is in $L^{\infty}(\mu)\bigcap\mathcal{A}_k^{\perp}$. Therefore for all $f \in L^{\infty}(\mu)$ we have
$$\frac{1}{N}\sum_{n=1}^N \left[f(T^{n}x)-\mathbb{E}(f|\mathcal{A}_k)(T^nx)\right]g_{1}(S_{1}^{n}y_{1})g_{2}(S_{2}^{n}y_{2})\cdots g_{k}(S_{k}^{n}y_{k})$$
converges to 0 $\nu_k$-a.e..  Thus, we have shown that the factors $\mathcal{A}_k$ are pointwise characteristic for the multiple term return times averages.

To finish the proof of the theorem it remains to show that
$$\limsup_N\left|\frac{1}{N}\sum_{n=1}^N f(T^nx)g_1(S_1^ny)g_{2}(S_{2}^{n}y_{2})\cdots g_{k}(S_{k}^{n}y_{k})\right|^2\leq CN_{k+1}(f)^2.$$
We use the property just established that the factors $\mathcal{A}_k$ are pointwise characteristic for the multiple term return times averages of $k+1$ functions including the arbitrary function $f$ and the Van der Corput lemma. We apply this characteristic property to each of the functions $f \cdot f\circ T^h$ and apply the Cauchy Schwartz inequality to obtain our result.

We have
\begin{eqnarray*}
\limsup_N\left|\frac{1}{N}\sum_{n=1}^N f(T^nx)g_1(S_1^ny)g_{2}(S_{2}^{n}y_{2})\cdots g_{k}(S_{k}^{n}y_{k})\right|^2 & \leq &\\
C \Bigg[ \frac{1}{H}+\frac{1}{H}\sum_{h=1}^H \limsup_N\Big| \frac{1}{N}\sum_{n=1}^{N-h}f(T^nx)f(T^{n+h}x)\cdot & & \\
g_{1}(S_{1}^{n}y_{1})g_{1}(S_{1}^{n+h}y_{1})\cdots g_{k}(S_{k}^{n}y_{k})g_{k}(S_{k}^{n+h}y_{k})\Big|\Bigg] & = &\\
C \Bigg[ \frac{1}{H} + \frac{1}{H}\sum_{h=1}^H \limsup_N\Big| \frac{1}{N}\sum_{n=1}^{N-h}\mathbb{E}(f \cdot f(T^h)|\mathcal{A}_k)(T^nx)\cdot & & \\ g_{1}(S_{1}^{n}y_{1})g_{1}(S_{1}^{n+h}y_{1})\cdots g_{k}(S_{k}^{n}y_{k})g_{k}(S_{k}^{n+h}y_{k})\Big|\Bigg] & &
\end{eqnarray*}

Applying the characteristic property to each of the functions $f \cdot f\circ T^h$  the above inequality is

\begin{equation}\label{Last2}
\leq   C\left[\frac{1}{H} + \frac{1}{H}\sum_{h=1}^H \left(\limsup_N\left|\frac{1}{N}\sum_{n=1}^{N-h}|\mathbb{E}(f \cdot f(T^h)|\mathcal{A}_k)(T^nx)|^2\right|\right)^{1/2}\right]
\end{equation}
because $\|g_i\|_{\infty}\leq 1$.

By Birkhoff's pointwise ergodic theorem and the ergodicity of $T$, the inequality in (\ref{Last2}) is
\begin{equation}\label{Last3}
\leq  C\left[\frac{1}{H} +\frac{1}{H}\sum_{h=1}^H \|\mathbb{E}(f\cdot f(T^h)|\mathcal{A}_k)\|\right].
\end{equation}

Using the Cauchy Schwartz Inequality we obtain that the inequality in (\ref{Last3}) is less than or equal to
$$C\left[\frac{1}{H} + \left(\frac{1}{H}\sum_{h=1}^H \|\mathbb{E}(f\cdot f(T^h)|\mathcal{A}_k)\|^2\right)^{1/2}\right].$$
Taking the limit with $H$ gives us the upper bound $C N_{k+1}(f)^2.$
\ep
\medbreak

\section{The $\mathcal{Z}_k$ factors are pointwise characteristic for the multiterm return times averages}

As noted above, the factors $\mathcal{Z}_k$ are smaller than the $\mathcal{A}_k$ factors and thus their orthogonal complements $\mathcal{Z}_k^{\perp}$ are bigger.  Therefore Theorem \ref{HKZ}, which we are proving in this section, is an extension of Theorem \ref{Main}.  We will prove Theorem  \ref{HKZ} directly from the properties of the factors $\mathcal{Z}_k$.  The Host-Kra-Ziegler factors, $\mathcal{Z}_k$, were defined in \cite{HK-NEA} by seminorms as follows.

\proc{Definition:}
\begin{itemize}
    \item The factor $\mathcal{Z}_0$ is equal to the trivial $\sigma$-algebra.
    \item The factor $\mathcal{Z}_1$ can be characterized by the seminorms $\||f|\|_2$  where
        $$\||f|\|_2^4 = \lim_{H}\frac{1}{H}\sum_{h=1}^{H} \left|\int f\cdot f\circ T^hd\mu\right|^2 $$
    \item The factor $\mathcal{Z}_2$ is the Conze-Lesigne factor, $\mathcal{CL}$.  Functions in this factor are characterized by the seminorm $|\|\cdot |\|_3$ such that
        $$\||f|\|_3^8 = \lim_{H}\frac{1}{H} \sum_{h=1}^{H}\||f\cdot f\circ T^h|\|_2^4.$$
        A function $f\in \mathcal{CL}^{\perp}$ if and only $\||f|\|_3 =0.$
    \item More generally B. Host and B. Kra showed in \cite{HK-NEA} that for each positive integer $k$ we have
        \begin{equation}\label{Semi}
            \||f|\|_{k+1}^{2^{k+1}}= \lim_H\frac{1}{H}\sum_{h=1}^{H}\||f\cdot f\circ T^h|\|_{k}^{2^k},
        \end{equation}
        with the condition that $f\in \mathcal{Z}_{k-1}^{\perp}$ if and only if $\||f|\|_{k}=0.$
\end{itemize}
\medbreak

One can compare the factors $\mathcal{Z}_k$ and $\mathcal{A}_k.$  First, the factors $\mathcal{A}_k$ are bigger than the factors $\mathcal{Z}_k$. More precisely we have the following.
\begin{itemize}
    \item The factors $\mathcal{A}_0$ and $\mathcal{Z}_0$ are equal to the trivial $\sigma$-algebra.
    \item The factors $\mathcal{A}_1$ and $\mathcal{Z}_1$ are also identical.  The seminorm $\||f|\|_2$ and $N_2(f)$ are equal.  Indeed
        $$\||f|\|_2^4 = \lim_{H}\frac{1}{H}\sum_{h=1}^{H} \left|\int f\cdot f\circ T^hd\mu\right|^2 = N_2(f)^4.$$
    \item The difference starts with the factors $\mathcal{A}_2$ and $\mathcal{Z}_2.$
        It is not difficult to find examples where $\mathcal{A}_2\neq \mathcal{Z}_2$ .  On the two torus the transformation $(x,y) \rightarrow (x+\alpha, y + \sqrt{\{x\}})$ where $\{x\}$ denotes the fractional part of $x,$ is an example for which the two factors differ.  More generally it can be shown that if the transformation on the two torus is given by $(x,y) \rightarrow (x+\alpha,y+\rho(x))$, where $\rho:\mathbb{T} \rightarrow \mathbb{T}$ is measurable, then $\mathcal{A}_2$ always coincides with the full algebra (\textit{i.e.}, the system is 2-step distal), and $\mathcal{Z}_2=\mathcal{A}_2$ only when $\rho$ is cohomologous to the affine co-cycle.
\end{itemize}

Note that the factors $\mathcal{Z}_k$ have a very rigid algebraic structure. They have the structure of a pro-nil system. See \cite{HK-NEA} for more details on the structure of these factors.

Our induction argument comes from reducing the return times averages by looking at an associated Wiener-Wintner type average using the following lemma.

\begin{lemma}\label{VDC2}
Let $(X, \mathcal{F}, \mu, T)$ be an ergodic dynamical system and $f\in L^{\infty}(\mu).$ Then for all positive integers $H$
we have
$$\limsup_N \sup_t \left|\frac{1}{N}\sum_{n=1}^N f(T^nx)e^{2\pi i nt}\right|^2 \leq C\left(\frac{1}{H} + \frac{1}{H}\sum_{h=1}^H \left|\int f\cdot f\circ T^hd\mu\right|\right)$$ where $C$ is an absolute constant derived from the application of the Van der Corput lemma.  In particular we have for $\mu$-a.e. $x$
$$\limsup_N \sup_t\left|\frac{1}{N}\sum_{n=1}^{N} f(T^nx) e^{2\pi int}\right|^2 \leq C \||f|\|_2^2.$$
\end{lemma}

\proc{Proof:}
This is Lemma 2 from the paper \cite{AssaniPCAC}.
\ep
\medbreak

Using this result, we can deduce the following lemma concerning the integral of the $\limsup$ of our averages.

\begin{lemma}\label{Induct}
Given $(X, \mathcal{F}, \mu, T)$ an ergodic measure preserving system on a probability measure space and $f\in L^{\infty}$.
Then we can find a set of full measure $X_f$ such that for every $x\in X_f$ for each measure preserving dynamical system $\Gamma_1=(Y_1, \mathcal{G}_1, \nu_1, S_1)$ and each $g_1\in L^{\infty}(\nu_1)$ with $\|g_1\|_{\infty} \leq 1$ we have
$$\int \limsup_N F_N^1(y_1) d\nu_1 \leq C \||f|\|_2^2$$
where $C$ is an absolute constant derived from the application of the Van der Corput lemma and
$$F_N^1(y_1):=\left|\frac{1}{N}\sum_{n=1}^{N}f(T^nx)g_1(S_1^ny_1)\right|^2.$$
\end{lemma}

\proc{Proof:}
By the BFKO \cite{RTT2} proof of the Return Times Theorem we have pointwise convergence of the above averages, therefore the $\limsup$ on the left hand side of the above expression becomes a limit.  Therefore, we have
\begin{eqnarray*}
\int \limsup_N F_N^1(y_1) d\nu_1 & = & \lim_N\int F_N^1(y_1)d\nu_1 \\
& = & \lim_N\int \left|\frac{1}{N}\sum_{n=1}^{N}f(T^nx)e^{2\pi int}\right|^2d\sigma_{g_1}(t)
\end{eqnarray*}
where $\sigma_{g_1}$ is the spectral measure associated to $g_1$ with respect to the dynamical system $\Gamma_1$.
Thus
$$\int \limsup_N F_N^1(y_1) d\nu_1 \leq \limsup_N\sup_t\left|\frac{1}{N}\sum_{n=1}^{N}f(T^nx)e^{2\pi int}\right|^2\|g_1\|_2^2.$$
As $\|g_1\|_{\infty} \leq 1$, using Lemma \ref{VDC2} we derive the inequality
$$\int \limsup_N F_N^1(y_1) d\nu_1 \leq C \||f|\|_2^2.$$
\ep
\medbreak

From Lemma \ref{Induct} the iteration process follows.  For instance, we can use this lemma to prove the following Wiener-Wintner return times result which refines the one obtained in \cite{ALR}.

\begin{lemma}\label{ALR}
Let $(X, \mathcal{F}, \mu, T)$ be an ergodic measure preserving system on a probability measure space and $f\in L^{\infty}(\mu)$.
Then for $\mu$-a.e. $x\in X$ for every measure preserving system $\Gamma_1 =(Y_1, \mathcal{G_1}, \nu_1, S_1)$ and each $g_1\in L^{\infty}(\nu_1)$ with $\|g_1\|_{\infty} \leq 1$ we have
\begin{equation}\label{LemmaALR}
\int \limsup_N \sup_t \left| \frac{1}{N} \sum_{n=1}^N f(T^nx)g_1(S_1^ny_1) e^{2\pi int}\right|^2 d\nu_1\leq C\||f|\|_3^2
\end{equation}
where $C$ is the absolute constant from the application of the Van der Corput lemma.  In particular, for $f\in \mathcal{CL}^{\perp}$ (or equivalently $\||f|\|_3 = 0$) we have for $\nu_1$-a.e. $y_1$
\begin{equation}\label{LemmaALR2}
\limsup_N \sup_t \left| \frac{1}{N} \sum_{n=1}^N f(T^nx)g_1(S_1^ny_1) e^{2\pi int}\right|= 0.
\end{equation}
\end{lemma}

\proc{Proof:}
By the Van der Corput lemma \cite{KN} we have
\begin{eqnarray}
\int \limsup_N \sup_t \left| \frac{1}{N} \sum_{n=1}^N f(T^nx)g_1(S_1^ny_1) e^{2\pi int}\right|^2d\nu_1 & \leq & \nonumber \\
C\Bigg(\frac{1}{H}+\frac{1}{H}\sum_{h=1}^H \int \limsup_N \bigg| \frac{1}{N}\sum_{n=1}^N f(T^nx)f(T^{n+h}x) \cdot & & \nonumber \\
g_1(S_1^ny_1)g_1(S_1^{n+h}y_1)\bigg|d\nu_1\Bigg) &  &\label{Ineq1}
\end{eqnarray}

Using the Cauchy-Schwarz inequality we have that the expression in (\ref{Ineq1}) is less than or equal to
\begin{eqnarray}\label{Ineq2}
C\Bigg(\frac{1}{H} + \frac{1}{H}\sum_{h=1}^H \bigg(\int \limsup_N \Big| \frac{1}{N}\sum_{n=1}^N f(T^nx)f(T^{n+h}x)\cdot \\ \nonumber
g_1(S_1^ny_1)g_1(S_1^{n+h}y_1)\Big|^2d\nu_1\bigg)^{1/2}\Bigg)
\end{eqnarray}

Similarly to the proof of Lemma \ref{Induct} as pointwise convergence of the averages holds by Theorem \ref{MRTT} we can rewrite the above $\limsup_N$ as a $\lim_N$ and use the spectral theorem to rewrite the integral in (\ref{Ineq2}) as
$$\lim_N\int \left|\frac{1}{N}\sum_{n=1}^{N}f(T^nx)f(T^{n+h}x)e^{2\pi int}\right|^2d\sigma_{g_1\cdot g_1\circ S_1^h}(t).$$
Thus by Lemma \ref{VDC2} the expression in \ref{Ineq2} is
\begin{eqnarray*}
& \leq & C\Bigg(\frac{1}{H}+\frac{1}{H}\sum_{h=1}^H \||f\cdot f\circ T^h|\|_2\|g_1\cdot g_1\circ S_1^h|\|_2\Bigg) \\
& \leq & C\Bigg(\frac{1}{H}+\frac{1}{H}\sum_{h=1}^H \||f\cdot f\circ T^h|\|_2\|g_1\|_{\infty}^2 \Bigg)\\
& \leq & C\Bigg(\frac{1}{H}+\frac{1}{H}\sum_{h=1}^H \||f\cdot f\circ T^h|\|_2\Bigg)
\end{eqnarray*}
on a set of full measure depending only on $f$ as $\|g\|_{\infty} \leq 1$.  This set of full measure is, in fact, the intersection of the sets of full measure obtained by the BFKO \cite{RTT2} proof of the Return Times Theorem for each function $f\cdot f\circ T^h.$

Continuing from above and using H\"{o}lder's Inequality we have
\begin{eqnarray*}
\int \limsup_N \sup_t \left| \frac{1}{N} \sum_{n=1}^N f(T^nx)g_1(S_1^ny_1) e^{2\pi int}\right|^2 d\nu_1 &\leq & \\
C\Bigg(\frac{1}{H}+\left(\frac{1}{H}\sum_{h=1}^H \||f\cdot f\circ T^h|\|_2^4\right)^{1/4}\Bigg). & &
\end{eqnarray*}

As the seminorm is defined by
$$\||f|\|_3^8 = \lim _H \frac{1}{H}\sum_{h=1}^H \||f\cdot f\circ T^h|\|_2^4$$
taking the limit on $H$ in the above expression gives
\begin{eqnarray*}
\int \limsup_N \sup_t \left| \frac{1}{N} \sum_{n=1}^N f(T^nx)g_1(S_1^ny_1) e^{2\pi int}\right|^2 d\nu_1 &\leq & C\left(\||f|\|_3^8\right)^{1/4} \\
& = &  C\||f|\|_3^2.
\end{eqnarray*}

This proves (\ref{LemmaALR}) of Lemma \ref{ALR}. Equation (\ref{LemmaALR2}) follows directly from the characterization of the $\mathcal{CL}$ factor.
\ep
\medbreak

The induction assumption giving the result on the pointwise characteristic factors for the $\mathcal{Z}_k$ factors can now be made.  To end it at the $\mathcal{CL}=\mathcal{Z}_2$ level we prove the next lemma.

\begin{lemma}\label{CL}
Let $(X, \mathcal{F}, \mu, T)$ be an ergodic measure preserving system on a probability measure space. The factor $\mathcal{Z}_2$, the Conze Lesigne factor, is pointwise characteristic  for the three term return times theorem.
\end{lemma}

\proc{Proof:}
We denote by $F_N(y_1,y_2)$ the three term averages with our original function $f$,  the fixed system $\Gamma_1=(Y_1, \mathcal{G}_1, \nu_1, S_1)$ and the variable one $\Gamma_2=(Y_2, \mathcal{G}_2, \nu_2, S_2).$ More precisely we have
$$F_N^2(y_1,y_2)=\left|\frac{1}{N}\sum_{n=1}^N f(T^nx)g_1(S_1^ny_1)g_2(S_2^ny_2)\right|^2$$
By Theorem \ref{MRTT} we have a set of full measure $Y_{g_1} \subset Y_1$ on which the pointwise convergence of the return times averages with three terms holds for any choice of measure preserving dynamical system $\Gamma_2=(Y_2, \mathcal{G}_2, \nu_2, S_2)$ and $g_2 \in L^{\infty}(\nu_2)$ with $\|g_2\|_{\infty} \leq 1$.  Therefore for $y_1 \in Y_{g_1}$ we have

\begin{equation}\label{lim}
\int\limsup_NF_N^2(y_1,y_2) d\nu_2=\lim_N \int F_N^2(y_1,y_2) d\nu_2
\end{equation}

Using the spectral measure as before, we can continue from (\ref{lim})
\begin{eqnarray*}
\int\limsup_NF_N^2(y_1,y_2) d\nu_2 & = & \lim_N \int \left|\frac{1}{N}\sum_{n=1}^N f(T^nx)g_1(S_1^ny_1)e^{2\pi int}\right|^2 d\sigma_{g_2}(t)\\
& \leq & \limsup_N \sup_t\left|\frac{1}{N}\sum_{n=1}^N f(T^nx)g_1(S_1^ny_1)e^{2\pi int}\right|^2\|g_2\|_{\infty}^2\\
& \leq &  \limsup_N \sup_t\left|\frac{1}{N}\sum_{n=1}^N f(T^nx)g_1(S_1^ny_1)e^{2\pi int}\right|^2
\end{eqnarray*}
as $\|g_2\|_{\infty}|\leq 1$.  Note that this upper bound is now independent of the choice of $\Gamma_2$ and $g_2$ so in fact we have
$$\sup_{\Gamma_2,g_2}\int \limsup_N F_N^2(y_1,y_2) d\nu_2 \leq \limsup_N \sup_t\left|\frac{1}{N}\sum_{n=1}^N f(T^nx)g_1(S_1^ny_1)e^{2\pi int}\right|^2.$$
Note that the left hand side of this last inequality is not necessarily measurable.  However one can conclude by making the following observation.
By Lemma \ref{ALR} we have
$$\int \limsup_N \sup_t \left| \frac{1}{N} \sum_{n=1}^N f(T^nx)g_1(S_1^ny_1) e^{2\pi int}\right|^2 d\nu_1\leq C\||f|\|_3^2.$$

Therefore, for $f \in \mathcal{Z}_2^{\perp}=\mathcal{CL}^{\perp}$ (i.e. $\||f|\|_3 = 0)$) there exists a set of full measure in $Y_1$ on which we have
$$\limsup_N \sup_t \left| \frac{1}{N} \sum_{n=1}^N f(T^nx)g_1(S_1^ny_1) e^{2\pi int}\right|=0.$$
For $y_1\in Y_1$ we have
$$\int \limsup_N \left|\frac{1}{N}\sum_{n=1}^N f(T^nx)g_1(S_1^ny_1)g_2(S_2^ny_2)\right|^2d\nu_2 = 0$$
for any choice of dynamical system $\Gamma_2$ and $g_2 \in L^{\infty}(\nu_2)$ with $\|g_2\|_{\infty} \leq 1$.  Therefore $\mathcal{Z}_2$ is pointwise characteristic for the three term return times averages.
\ep
\medbreak

We now have the tools necessary to prove our second main result, Theorem \ref{HKZ}, that the $\mathcal{Z}_k$ averages are pointwise characteristic for the multiterm return times averages.

\proc{Proof:}
It remains to finish the induction argument which we have started in the above lemmas.  Suppose that we know that the $\mathcal{Z}_{j}$ are pointwise characteristic for the $j$-th return times averages and that the bound in (\ref{LemmaALR}) from Lemma \ref{ALR} holds for $j$-th averages for $1 \leq j <k$.

Specifically for $\mu$-a.e. $x\in X$, for every measure preserving system $\Gamma_j=(Y_j, \mathcal{G}_j, \nu_j, S_j)$ and each $g_j\in L^{\infty}(\nu_j)$ with $\|g_j\|_{\infty} \leq 1$ with $1 \leq j < k-1$ we have for $f\in \mathcal{Z}_k^{\perp}$
\begin{equation}\label{CLIneq2}
  \limsup_N \Delta_N^{k-1}(y_1,\ldots,y_{k-1}) = 0
\end{equation}
where
$$\Delta_N^{k-1}(y_1,\ldots, y_{k-1})=\sup_t\left|\frac{1}{N}\sum_{n=1}^N f(T^nx)g_1(S_1^ny_1)\cdots g_{k-1}(S_{k-1}^ny_{k-1})e^{2\pi int}\right|^2$$
and $\Gamma_{k-1}=(Y_{k-1}, \mathcal{G}_{k-1}, \nu_{k-1}, S_{k-1})$ is any measure preserving dynamical system and $g_{k-1} \in L^{\infty}(\nu_{k-1})$ with $\|g_{k-1}\|_{\infty} \leq 1$.

By Theorem \ref{MRTT}, we know that for any $f \in L^{\infty}(\mu)$ and $1 \leq j < k$ we can find sets $X_f$ of full measure such that if $x \in X_f$, then for any other dynamical system $\Gamma_1=(Y_1,\mathcal{G}_1, S_1, \nu_1)$ and any $g_1 \in L^{\infty}(\nu_1)$ with $\|g_1\|_{\infty} \leq 1$, there exists a set of full measure $Y_{g_1}$ such that for each $y_1$ in $Y_{g_1}$ then $\ldots$ for any other dynamical system $\Gamma_{k-1}=(Y_{k-1},\mathcal{G}_{k-1}, S_{k-1}, \nu_{k-1})$ and any $g_{k-1} \in L^{\infty}(\nu_{k-1})$ with $\|g_{k-1}\|_{\infty} \leq 1$ there exist a set of full measure $Y_{g_{k-1}}$ in $Y_{k-1}$ such that if $y_{k-1} \in Y_{g_{k-1}}$ we have the pointwise convergence of the return times averages with $k$ terms for any other dynamical system $\Gamma_k=(Y_k,\mathcal{G}_k, S_k, \nu_k)$ and any $g_{j} \in L^{\infty}(\nu_{j})$ with $\|g_{j}\|_{\infty} \leq 1$.

Thus for $x,y_1,\ldots y_{k-1}$ as above we have
\begin{equation}\label{limk}
\int\limsup_NF_N^{k}(y_1,\ldots,y_{k}) d\nu_k=\lim_N \int F_N^{k}(y_1,\ldots,y_{k}) d\nu_k
\end{equation}
where
$$F_N^k(y_1,\ldots,y_{k})=\left|\frac{1}{N}\sum_{n=1}^N f(T^nx)g_1(S_1^ny_1)...g_{k}(S_{k}^ny_{k})\right|^2.$$

We would like to show that $\mathcal{Z}_{k}$ is pointwise characteristic for the $k$-return times averages and that we have for $f \in \mathcal{Z}_{k+1}^{\perp}$
$$
  \limsup_N \Delta_N^{k}(y_1,\ldots,y_{k})=0
$$
where the
$$\Delta_N^{k}(y_1,\ldots, y_{k})=\sup_t\left|\frac{1}{N}\sum_{n=1}^N f(T^nx)g_1(S_1^ny_1)\cdots g_{k}(S_{k}^ny_{k})e^{2\pi int}\right|^2,$$
$\Gamma_k$ and $g_k$ are as defined above.  This will complete our proof by induction.

Using the spectral theorem and continuing from (\ref{limk}) we have
\begin{eqnarray*}
\int\limsup_N F_N^{k}(y_1,\ldots,y_{k}) d\nu_k & = & \\
\lim_N\int \bigg|\frac{1}{N}\sum_{n=1}^N f(T^nx)g_1(S_1^ny_1)\cdots g_{k-1}(S_{k-1}^ny_{k-1})e^{2\pi int}\bigg|^2d\sigma_{g_k}(t) & \leq & \\
\limsup_N\sup_t\bigg|\frac{1}{N}\sum_{n=1}^N f(T^nx)g_1(S_1^ny_1)\cdots g_{k-1}(S_{k-1}^ny_{k-1})e^{2\pi int}\bigg|^2\|g_k\|_{\infty} & \leq & \\
\limsup_N\sup_t\bigg|\frac{1}{N}\sum_{n=1}^N f(T^nx)g_1(S_1^ny_1)\cdots g_{k-1}(S_{k-1}^ny_{k-1})e^{2\pi int}\bigg|^2 &  &
\end{eqnarray*}
as $\|g_k\|\leq 1$.

Note that this upper bound is now independent of the choice of $\Gamma_k$ and $g_k$ so
$$\sup_{\Gamma_k,g_k}\int\limsup_NF_N^{k}(y_1,\ldots,y_{k}) d\nu_k \leq $$
\begin{equation}\label{Ugly}
\limsup_N\sup_t\bigg|\frac{1}{N}\sum_{n=1}^N f(T^nx)g_1(S_1^ny_1)\cdots g_{k-1}(S_{k-1}^ny_{k-1})e^{2\pi int}\bigg|^2.
\end{equation}

In the same manner as shown in Lemma \ref{CL}, from the equation (\ref{CLIneq2}) and (\ref{Ugly}) one can conclude that $\mathcal{Z}_k$ is pointwise characteristic for the $k$-th return times averages.

To complete the induction step it remains to show that (\ref{CLIneq2}) holds for $k$.

To make the reading of the induction proof easier we show how one can prove that $\mathcal{Z}_4$ is pointwise characteristic for the $5$-th return times averages. The reader will check that the arguments extend without difficulty to arbitrary $k$.
So we want to show that if $\||f|\|_5 =0$ then
\begin{equation}\label{Ind5}
\limsup_N\sup_t\big|\frac{1}{N}\sum_{n=1}^N f(T^nx)g_1(S_1^ny_1)\cdots g_4(S_4^ny_4)e^{2\pi int}\big|^2 =0.
\end{equation}

By the Van der Corput lemma, Theorem \ref{MRTT} and Cauchy-Schwartz inequality we have
\begin{eqnarray*}
& & \int\limsup_N\sup_t\left|\frac{1}{N}\sum_{n=1}^N f(T^nx)g_1(S_1^ny_1)\cdots g_4(S_4^ny_4)e^{2\pi int}\right|^2d\nu_4 \\
& \leq & C\Bigg(\frac{1}{H_1} + \Big(\frac{1}{H_1}\sum_{h_1=1}^{H_1} \limsup_N\int \big|\frac{1}{N}\sum_{n=1}^N f(T^nx)f(T^{n+h_1}x)\cdot \\
& & \cdots g_4(S_4^ny_4)g_4(S_4^{n+h_1}y_4)\big|^2d\nu_4\Big)^{1/2}\Bigg).
\end{eqnarray*}

By the spectral theorem this last term is equal to
\begin{eqnarray*}
& & \Big(\frac{1}{H_1}\sum_{h_1=1}^{H_1} \limsup_N\int \big|\frac{1}{N}\sum_{n=1}^N f(T^nx)f(T^{n+h_1}x)\cdot \\
& & \cdots g_3(S_3^ny_3)g_3(S_3^{n+h_1}y_3)e^{2\pi int}\big|^2d\sigma_{g_4.g_4\circ S_4^{h_1}}\Big)^{1/2}.
\end{eqnarray*}

As $\|g_4\|_{\infty}\leq 1,$ the generic term in this sum is less than
$$\limsup_N\sup_t\big|\frac{1}{N}\sum_{n=1}^N f(T^nx)f(T^{n+h_1}x)\cdots g_3(S_3^ny_3)g_3(S_3^{n+h_1}y_3)e^{2\pi int}\big|^2.$$

One can conclude that for the appropriate universal sets for $x,$ $y_1,$ $y_2,$ and $y_3$ if
\begin{eqnarray*}
\lim_{H_1}\frac{1}{H_1}\sum_{h_1=1}^{H_1} \limsup_N\sup_t\big|\frac{1}{N}\sum_{n=1}^N f(T^nx)f(T^{n+h_1}x)\cdot & &  \\
\cdots g_3(S_3^ny_3)g_3(S_3^{n+h_1}y_3)e^{2\pi int}\big|^2 & = & 0
\end{eqnarray*}
then
$$\int\limsup_N\sup_t\big|\frac{1}{N}\sum_{n=1}^N f(T^nx)g_1(S_1^ny_1)\cdots g_4(S_4^ny_4)e^{2\pi int}\big|^2d\nu_4=0$$
for all measure preserving system $\Gamma_4=(Y_4, \mathcal{G}_4, \nu_4, S_4)$ and each $g_4$ with $\|g_4|\|_{\infty}\leq 1.$

We then derive that
$$\limsup_N\sup_t\big|\frac{1}{N}\sum_{n=1}^N f(T^nx)g_1(S_1^ny_1)\cdots g_4(S_4^ny_4)e^{2\pi int}\big|^2=0$$
which is the content of (\ref{Ind5}).

Therefore we look at the quantity
$$\int\limsup_N\sup_t\big|\frac{1}{N}\sum_{n=1}^N f(T^nx)f(T^{n+h_1}x)\cdots g_3(S_3^ny_3)g_3(S_3^{n+h_1}y_3)e^{2\pi int}\big|^2d\nu_3.$$
Again by Van der Corput lemma, Theorem \ref{MRTT} and Cauchy Schwartz inequality this quantity is less than
\begin{eqnarray*}
C\Bigg(\frac{1}{H_2} + \Big(\frac{1}{H_2}\sum_{h_2=1}^{H_2} \limsup_N\int \big|\frac{1}{N}\sum_{n=1}^N f(T^nx)f(T^{n+h_1}x)f(T^{n+h_2}x) f(T^{n+h_1+h_2}x) \\ \cdots g_3(S_3^ny_3)g_3(S_3^{n+h_1}y_3)g_3(S_3^{n+h_2}y_3)g_3(S_3^{n+h_1+h_2}y_3)\big|^2d\nu_3\Big)^{1/2}\Bigg). \end{eqnarray*}

As shown previously, by the spectral theorem this last term is equal to
\begin{eqnarray*}
\Big(\frac{1}{H_2}\sum_{h_2=1}^{H_2} \limsup_N\int \big|\frac{1}{N}\sum_{n=1}^N f(T^nx)f(T^{n+h_1}x)f(T^{n+h_2}x) f(T^{n+h_1+h_2}x)\cdot \\
g_2(S_2^ny_2)g_2(S_2^{n+h_1}y_2)g_2(S_2^{n+h_2}y_2)g_2(S_2^{n+h_1+h_2}y_2)\cdot \\ e^{2\pi int}\big|^2d\sigma_{g_3\cdot g_3\circ S_3^{h_1}g_3\circ S_3^{h_2}g_3\circ S_3^{h_1+h_2}}\Big)^{1/2}.
\end{eqnarray*}

This in turn is less than
\begin{eqnarray*}
\frac{1}{H_2}\sum_{h_2=1}^{H_2} \limsup_N\sup_t \big|\frac{1}{N}\sum_{n=1}^N f(T^nx)f(T^{n+h_1}x)f(T^{n+h_2}x) f(T^{n+h_1+h_2}x)\cdot \\
g_2(S_2^ny_2)g_2(S_2^{n+h_1}y_2)g_2(S_2^{n+h_2}y_2)g_2(S_2^{n+h_1+h_2}y_2)e^{2\pi int}\big|^2.
\end{eqnarray*}

We integrate this term with respect to $\nu_2$. Another application of the Van der Corput lemma, Theorem \ref{MRTT}, Cauchy Schwartz Inequality and the spectral theorem leads to the estimate

$$\Big(\frac{1}{H_3}\sum_{h_3=1}^{H_3} \limsup_N\sup_t\big|\frac{1}{N}\sum_{n=1}^N f(T^nx)f(T^{n+h_1}x)f(T^{n+h_2}x) f(T^{n+h_1+h_2}x)$$
$$f(T^{n+h_3}x)f(T^{n+h_1+h_3}x)f(T^{n+h_2+h_3}x)f(T^{n+h_1+h_2+h_3}x) g_1(S_1^ny_1)g_1(S_1^{n+h_1}y_1)$$
$$g_1(S_1^{n+h_2}y_1)g_1(S_1^{n+h_1+h_2}y_1)g_1(S_1^{n+h_3}y_1)g_1(S_1^{n+h_1+h_3}y_1)g_1(S_1^{n+h_2+h_3}y_1)$$
$$g_1(S_1^{n+h_1+h_2+h_3}y_1)e^{2\pi int}\big|^2\Big)^{1/2}.$$

Finally integrating this last term with respect to $\nu_1$ we can use (\ref{LemmaALR}) and combine the previous inequalities to get the upper bound
\begin{eqnarray*}
C\Bigg(\sum_{i=1}^3\frac{1}{H_i} +
\frac{1}{\prod_{i=1}^3 H_i}\sum_{h_3=1}^{H_3}\sum_{h_2=1}^{H_2}\sum_{h_1=1}^{H_1}
\||f\cdot f\circ T^{h_1}\cdot & & \\
f\circ T^{h_2}f\circ T^{h_1+h_2}\cdot f\circ T^{h_3}f\circ T^{h_1+h_3}f\circ T^{h_2+h_3}f\circ T^{h_1+h_2+h_3}|\|_3\Bigg)& &
\end{eqnarray*}

Applying (\ref{Semi}) three times gives the upper bound
$$ C\||f|\|_{5}^8.$$
Thus if $\||f|\|_{5}=0$ we can go back and step by step obtain universal sets for $x,y_1, y_2, y_3$  for which (\ref{Ind5}) holds.
\ep
\medbreak

\proc{Remarks:}
\begin{enumerate}
\item In Theorem \ref{A_k} (see Equation (\ref{Upbound})), we proved that we have a pointwise upper bound on the average of multiple terms as follows
$$\limsup_N \left| \frac{1}{N}\sum_{n=1}^N f(T^{n}x)g_{1}(S_{1}^{n}y_{1})\cdots g_{k}(S_{k}^{n}y_{k})\right|^2 \leq C N_{k+1}(f)^2.$$  We asked in a previous version of this paper whether for $k\geq 2,$ one can replace in these inequalities the $N_k$ seminorms for the $\mathcal{A}_k$ with those defining the $\mathcal{Z}_k$ factors.
\begin{itemize}
\item It was shown in \cite{HK-NEA} that $\||f|\|_2\leq \|\mathbb{E}(f|\mathcal{K}_T)\|_2$ but we can not find an absolute constant $C$ for which $\|\mathbb{E}(f|\mathcal{K}_T)\|_2\leq C \||f|\|_2.$ 
\item As pointed out by the referee one can not replace the $N_k$ seminorms in Equation (\ref{Upbound}) with those defining the $\mathcal{Z}_k$ factors. For $T$ an irrational rotation on the circle the Kronecker factor is $\mathcal{F}$ and $\mathbb{E}(f, \mathcal{K}_T) = f$. Then $\|\mathbb{E}(f| \mathcal{K}_{T})\|_2 = \|\hat{f}(k)\|_{l^2\mathbb(Z)}$. Direct computations show that $\||f|\|_2= \|\hat{f}(k)\|_{l^4\mathbb(Z)}.$ Therefore one can not find an absolute constant $C$ for which
    $\|\hat{f}\|_{l^4(\mathbb{Z})}\leq C\|\hat{f}\|_{l^2(\mathbb{Z})}.$ Consider the average $\frac{1}{N}\sum_{n=1}^N f(T^nx)\overline{f}(T^ny).$ Then we have 
    $$\lim_N \left|\frac{1}{N}\sum_{n=1}^N f(T^nx)\overline{f}(T^ny)\right| = \left|\sum_{k\in \mathbb{Z}}|\hat{f}(k)|^2e^{2\pi ik(x-y)}\right|.$$ 
    The right hand side of the above inequality is less than $\sum_{k\in\mathbb{Z}}|\hat{f}(k)|^2= \|\mathbb{E}(f| \mathcal{K}_{T})\|_2^2.$ However one can not find an absolute constant $C$ for which
     $\big|\sum_{k\in \mathbb{Z}}|\hat{f}(k)|^2e^{2\pi ik(x-y)}\big|\le C \|\hat{f}(k)\|_{l^4\mathbb(Z)}^2.$
\end{itemize}
\item The authors of this paper are writing a survey of the Return Times Theorem \cite{AP} which will include more details of the historical developments of Theorem \ref{RTT1} and \ref{MRTT} and related questions such as the ones noted above.
\end{enumerate}
\noindent{\bf Acknowledgements} The authors thank the referees for their comments.
\medbreak

\end{document}